\ifdefined\pdfoutput\pdfoutput=1\fi
\documentclass[11pt,a4paper]{article}

\ifdefined\luatexversion\else\usepackage[utf8]{inputenc}\fi
\usepackage[T1]{fontenc}
\usepackage[english]{babel}
\usepackage{lmodern}
\usepackage{microtype}

\usepackage{amsmath,amssymb,amsfonts,amsthm}
\usepackage{geometry}
\usepackage{booktabs}
\usepackage{array}
\usepackage{caption}
\usepackage{cite}
\usepackage{graphicx}

\usepackage{hyperref}
\hypersetup{
    pdftitle={A Method for Linearization of Systems of Quadratic Equations: Local Convergence and Differentiable Physics Implicit Layers},
    pdfauthor={A. A. Pak, A. S. Pak},
    pdfkeywords={Linearization, Quadratic Systems, Kantorovich Inequality, Nonlinear Preconditioning, Deep Equilibrium Models, Differentiable Physics, Scientific Machine Learning},
    colorlinks=true,
    linkcolor=blue,
    citecolor=blue,
    urlcolor=blue
}

\newtheorem{theorem}{Theorem}[section]

\newtheorem{proposition}[theorem]{Proposition}
\newtheorem{corollary}[theorem]{Corollary}
\theoremstyle{definition}
\newtheorem{definition}[theorem]{Definition}
\newtheorem{remark}[theorem]{Remark}

\numberwithin{equation}{section}

\begin{document}

\title{\textbf{A Method for Linearization of Systems of Quadratic Equations:\\[0.6ex]
\large Theoretical Foundations, Local Convergence, and Differentiable Physics Implicit Layers}}

\author{
\textbf{A.~A.~Pak}\thanks{Corresponding author. E-mail: [\texttt{andr3y.pak@gmail.com}]}
\quad and \quad
\textbf{A.~S.~Pak}
}

\date{}

\maketitle

\begin{abstract}
Integrating physical conservation laws into deep neural architectures is a central challenge in Scientific Machine Learning (SciML). Soft-penalty formulations, including many Physics-Informed Neural Networks (PINNs), can be sensitive to loss weighting and optimization choices; generic Deep Equilibrium Model (DEQ) and implicit-layer solvers can likewise have instance-dependent iteration counts. This paper studies a norm-linearization scheme for networked systems governed by linear mass continuity and nonlinear quadratic or power-law energy/head-loss balances.

We first derive the sharp universal relative-variance bound $\eta \leq \frac{(p-1)^2}{4p}$ (alongside its exact finite-$n$ discrete counterpart) and connect it directly, via the weighted Kantorovich inequality, to the worst-case single-step squared-norm distortion of the linearization operator when the target state lies within a multiplicative condition ratio $p$ of the operating point. By applying Euler's theorem for degree-one homogeneous functions to directional norm balances, we then show that the linearization scheme is algebraically identical to a Newton--Raphson step applied to a system of norm differences $u_i(\mathbf{x}) - v_i(\mathbf{x}) = 0$. Rewriting raw degree-$2$ quadratic (or degree-$\alpha$ power-law) cycle balances $u_i(\mathbf{x})^\alpha - v_i(\mathbf{x})^\alpha = 0$ as degree-one homogeneous norm differences acts as a nonlinear degree-reducing preconditioner: it eliminates radial curvature ($\nabla^2 F_i(\mathbf{x})\mathbf{x} = \mathbf{0}$), renders the Jacobian scale-invariant (degree-zero homogeneous), causes the affine constant term in the Newton step to vanish identically, and automatically row-equilibrates the loop Jacobian at equilibrium. We prove that the coupled continuity-loop Jacobian is nonsingular at every regular equilibrium, establishing local quadratic convergence across both Darcy--Weisbach ($\alpha=2$) and Hazen--Williams ($\alpha=1.852$) regimes. Finally, we embed this solver with dynamic sign partitioning into a Differentiable Physics Implicit Layer (\texttt{FlowDEQ}) using a reduced incidence matrix and closed-form Implicit Function Theorem (IFT) adjoint differentiation.

\vspace{1ex}
\noindent\textbf{Keywords:} Linearization; Quadratic and Power-Law Systems; Kantorovich Inequality; Nonlinear Preconditioning; Deep Equilibrium Models; Differentiable Physics; Scientific Machine Learning.

\vspace{0.5ex}
\noindent\textbf{MSC 2020:} 65H10, 90C35, 68T07, 49M15.
\end{abstract}

\section{Introduction}

The integration of governing physical laws with modern deep learning architectures has become one of the most vital frontiers in computational engineering and artificial intelligence \cite{raissi2019,bai2019,amos2017,chen2018}. In physical systems defined on graphs and spatial networks, such as municipal fluid distribution grids \cite{todini1987}, regional natural gas pipelines, power transmission networks, and arterial vascular hemodynamics in cardiovascular medicine, the physical state is governed by two fundamentally distinct principles:
\begin{enumerate}
\item \textbf{Linear Conservation at Junction Nodes:} At each graph node $j \in \mathcal{V}$, the algebraic sum of branch flows must balance external injections or consumer demands:
\begin{equation}
\mathbf{A} \mathbf{x} = \mathbf{q}, \quad \text{i.e.}, \quad \sum_{k \in \operatorname{adj}(j)} \pm x_k - q_j = 0, \label{eq:continuity_intro}
\end{equation}
where $\mathbf{A} \in \mathbb{R}^{m \times n}$ is the incidence matrix, $\mathbf{x} \in \mathbb{R}^n$ represents branch flow rates, and $\mathbf{q} \in \mathbb{R}^m$ denotes nodal boundary conditions.
\item \textbf{Nonlinear Conservation around Closed Circuits:} For each closed cycle (loop) $l \in \mathcal{L}$ in the graph, the cumulative potential or frictional dissipation must sum to zero:
\begin{equation}
\sum_{k \in \operatorname{loop}(l)} s_k |x_k|^{\alpha-1} x_k = 0, \label{eq:energy_intro}
\end{equation}
where $s_k > 0$ represents physical dissipation/resistance coefficients, and $\alpha \ge 1$ is the flow friction exponent ($\alpha = 2$ for the classical quadratic Darcy--Weisbach law and high-Reynolds-number regimes; $\alpha \approx 1.852$ for the Hazen--Williams empirical relation).
\end{enumerate}

\subsection{Limitations of Current Machine Learning Paradigms}
When neural networks are tasked with solving inverse problems, parameter estimation, or digital-twin state estimation on such networks, two primary deep learning paradigms are currently deployed, both of which suffer from structural shortcomings:

\subsubsection{1. Soft-Penalty Physics-Informed Neural Networks (PINNs)}
In the classical PINN framework \cite{raissi2019}, a neural network parameterizes the flow state $\mathbf{x}_{\mathbf{w}}(\mathbf{q}_r)$ and physical laws are enforced as soft penalty terms in the scalar loss:
\begin{equation}
\mathcal{L}(\mathbf{w}, \mathbf{z}) = \mathcal{L}_{\text{data}}(\mathbf{x}_{\mathbf{w}}) + \lambda_{\text{node}} \|\mathbf{A}_r \mathbf{x}_{\mathbf{w}} - \mathbf{q}_r\|_2^2 + \lambda_{\text{loop}} \sum_{l \in \mathcal{L}} \left( \sum_{k \in \operatorname{loop}(l)} s_k(\mathbf{z}) |x_k|^{\alpha-1} x_k \right)^2.
\end{equation}
Despite their simplicity, soft-penalty PINNs can encounter optimization and scaling difficulties on cyclic networks:
\begin{itemize}
\item The loss landscape can be stiff, multi-scale, and non-convex.
\item Data fidelity and conservation penalties can have competing scales and gradient directions.
\item Constraint accuracy depends on architecture capacity, scaling, penalty weights, optimizer, and training budget; poorly calibrated penalties leave material conservation residuals.
\end{itemize}

\subsubsection{2. Generic Deep Equilibrium Models and Implicit Layers}
To guarantee exact constraint satisfaction, recent advances in \textbf{Deep Equilibrium Models (DEQs)} \cite{bai2019}, continuous-depth implicit models \cite{chen2018}, and \textbf{differentiable optimization layers (OptNet)} \cite{amos2017} formulate the physical state as an implicit root:
\begin{equation}
\mathbf{x}^* = \operatorname{Solve}(\mathbf{F}(\mathbf{x}^*, \boldsymbol{\theta}) = \mathbf{0}).
\end{equation}
Generic DEQ implementations may rely on off-the-shelf root-finding solvers (e.g., Broyden's quasi-Newton method or standard Newton applied directly to polynomial conservation equations \cite{todini1987}). On coupled quadratic and power-law network manifolds:
\begin{itemize}
\item Iteration counts can be instance-dependent, reducing utilization in batched GPU implementations.
\item Standard Newton applied directly to raw degree-$2$ quadratic (or degree-$\alpha$ power-law) loop equations suffers from radial curvature and non-vanishing affine shifts far from the root, failing to converge under tight fixed-step budgets or requiring line-search damping.
\end{itemize}

\subsection{Contributions of this Paper}
In this paper, we address both paradigms by presenting a tailored \textbf{Differentiable Physics Implicit Layer} based on an analytical norm-linearization of systems of quadratic and power-law equations.

The core contributions of this work are organized as follows:
\begin{enumerate}
\item \textbf{Extremum Analysis and Kantorovich Distortion Bounds:} We derive the sharp universal upper bound $\eta \leq \frac{(p-1)^2}{4p}$ for the relative variance functional, state the exact finite-$n$ equal-weight maximum, and prove via the weighted Kantorovich inequality that this bound governs the exact worst-case single-step squared-norm approximation distortion of our linearization scheme.
\item \textbf{Analytical Norm-Linearization Scheme:} We develop a linearization technique that transforms coupled quadratic and power-law loop equations into an unshifted, globally coupled system of linear equations solved simultaneously across all circuits at each iteration with dynamic sign partitioning $\boldsymbol{\sigma}^{(k)} = \operatorname{sgn}(\mathbf{x}^{(k)})$.
\item \textbf{Newton Equivalence, Nonsingularity, and Nonlinear Preconditioning:} Using Euler's homogeneous-function theorem, we prove that the linearization is algebraically identical to Newton--Raphson applied to degree-one homogeneous norm differences $u_i(\mathbf{x}) - v_i(\mathbf{x}) = 0$. We prove that the coupled Jacobian is nonsingular at every regular equilibrium via an exact row-equilibration identity $\mathbf{J}_{\text{loop}}(\mathbf{x}^*) = \mathbf{W}^{-(\alpha-1)}\mathbf{C}\mathbf{D}(\mathbf{x}^*)$, and show that transforming $u_i^\alpha - v_i^\alpha = 0 \mapsto u_i - v_i = 0$ acts as a nonlinear degree-reducing preconditioner that eliminates radial curvature ($\nabla^2 F_i(\mathbf{x})\mathbf{x} = \mathbf{0}$), normalizes Jacobian scaling, and eliminates the affine constant term.
\item \textbf{Power-Law Extension and Multi-Scale Verification:} We generalize the norm-difference formulation to weighted $\ell_\alpha$ norms for $\alpha>1$ (including Hazen--Williams $\alpha=1.852$) and verify that $K=4$ fixed-budget convergence holds across both power-law exponents and larger multi-loop networks ($n=113$ branches, $c=50$ loops).
\item \textbf{Differentiable Physics Layer (\texttt{FlowDEQ}):} We prove global existence and uniqueness of the equilibrium state via a strictly convex coercive variational principle, derive closed-form Implicit Function Theorem (IFT) adjoint gradients, and demonstrate robust inverse parameter recovery under both noiseless and noisy sparse sensor observations.
\end{enumerate}

\section{Extremum Analysis and the Error Bound \texorpdfstring{$\eta_{\max}$}{eta max}}

\subsection{Problem Formulation and Normalization}

Let $x_1, x_2, \dots, x_n$ be positive quantities ($x_i > 0$). Without loss of generality, we order them monotonically \cite{bronshtein1981}:
\begin{equation}
x_1 \le x_2 \le x_3 \le \dots \le x_n.
\end{equation}
If the elements are not ordered, they are re-indexed to satisfy this condition. Normalizing with respect to the minimum element $x_1$, we define:
\begin{equation}
\mu_i = \frac{x_i}{x_1}, \quad i = 1, 2, \dots, n,
\end{equation}
which implies:
\begin{equation}
\mu_1 = 1 \le \mu_2 \le \mu_3 \le \dots \le \mu_n = p, \quad \text{where } p = \frac{x_n}{x_1} \ge 1.
\end{equation}

Consider the relative variance functional:
\begin{equation}
\eta = \frac{n \sum_{i=1}^n x_i^2 - \left(\sum_{i=1}^n x_i\right)^2}{\left(\sum_{i=1}^n x_i\right)^2} = \frac{n \sum_{i=1}^n \mu_i^2}{\left(\sum_{i=1}^n \mu_i\right)^2} - 1.
\end{equation}
To determine the maximum distortion $\eta_{\max}$, it suffices to maximize the auxiliary objective:
\begin{equation}
I(\mu_1, \dots, \mu_n) = \frac{n \sum_{i=1}^n \mu_i^2}{\left(\sum_{i=1}^n \mu_i\right)^2}, \quad \text{such that } \eta = I - 1,
\end{equation}
subject to the box constraints \cite{vasiliev1980}:
\begin{equation}
1 \le \mu_i \le p, \quad i = 1, 2, \dots, n.
\end{equation}

\subsection{Reduction to Boundary Values}

For $n\ge2$, fix all coordinates except $\mu_i=t$, and set
\begin{equation}
B=\sum_{j\ne i}\mu_j>0, \qquad Q=\sum_{j\ne i}\mu_j^2>0.
\end{equation}
Apart from the positive constant $n$, the resulting one-variable objective is
\begin{equation}
g(t)=\frac{Q+t^2}{(B+t)^2}, \qquad 1\le t\le p,
\end{equation}
with derivative
\begin{equation}
g'(t)=\frac{2(Bt-Q)}{(B+t)^3}.
\end{equation}
Thus its only interior critical point, when it lies in the interval, is a minimum rather than a maximum. Holding the other coordinates fixed, a maximum is therefore attained at $t=1$ or $t=p$. Repeating this coordinatewise shows that a global maximum of $I$ on $[1,p]^n$ is attained at a vertex. (For $n=1$, $\eta\equiv0$.) Hence candidate maximizers have the form
\begin{equation}
\boldsymbol{\mu}^* = (\underbrace{1, 1, \dots, 1}_{n-k}, \underbrace{p, p, \dots, p}_{k}), \quad k \in \{1, 2, \dots, n-1\}.
\end{equation}

\subsection{Analytical Derivation of the Maximum}

For an extreme point with $n-k$ elements equal to $1$ and $k$ elements equal to $p$:
\begin{align}
\sum_{i=1}^n \mu_i &= (n - k) \cdot 1 + k \cdot p = n + k(p - 1), \\
\sum_{i=1}^n \mu_i^2 &= (n - k) \cdot 1^2 + k \cdot p^2 = n + k(p^2 - 1).
\end{align}
The auxiliary functional evaluated on this family of vertices is:
\begin{equation}
I(n, k, p) = \frac{n [n + k(p^2 - 1)]}{[n + k(p - 1)]^2}. \label{eq:I_discrete}
\end{equation}
Treating $k$ as a continuous relaxation parameter in the interval $[0, n]$:
\begin{equation}
\frac{dI}{dk} = \frac{n(p-1)^2 [n-k(p+1)]}{[n + k(p - 1)]^3}. \label{eq:I_derivative}
\end{equation}
Setting the derivative to zero yields the unique critical point:
\begin{equation}
\widehat{k} = \frac{n}{p+1}.
\end{equation}
For $p>1$, equation \eqref{eq:I_derivative} is positive for $k<\widehat{k}$ and negative for $k>\widehat{k}$, so this is the continuous maximum. Substitution gives the sharp universal bound
\begin{equation}
I \leq \frac{(p+1)^2}{4p}. \label{eq:I_univ}
\end{equation}
For a fixed finite $n$ with equal weights, $k$ must be an integer. The exact discrete maximum is instead
\begin{equation}
I_{n,\max}(p)=\max\!\left\{I\!\left(n,\left\lfloor\widehat{k}\right\rfloor,p\right), I\!\left(n,\left\lceil\widehat{k}\right\rceil,p\right)\right\},
\end{equation}
where repeated candidates are counted once. Consequently $I_{n,\max}(p)\le (p+1)^2/(4p)$, with equality precisely when $n/(p+1)$ is an integer (apart from the trivial $p=1$ case). Thus the universal relative-variance bound is
\begin{equation}
\eta \leq \eta_{\mathrm{univ}}(p) \equiv \frac{(p-1)^2}{4p}, \label{eq:etamax}
\end{equation}
whereas the exact finite-$n$ value is $\eta_{n,\max}(p)=I_{n,\max}(p)-1$.

\subsection{Connection to the Kantorovich and P\'olya--Szeg\H{o} Inequalities}

\begin{remark}
\label{rem:kantorovich}
The upper bound \eqref{eq:etamax} connects directly to two classical inequalities in functional analysis and provides the exact weighted generalization needed to analyze the quadratic linearization scheme in Section~3:
\begin{enumerate}
\item \textbf{The P\'olya--Szeg\H{o} Inequality (1925):} For two sets of positive numbers $a_i, b_i$ with $0 < m_1 \le a_i \le M_1$ and $0 < m_2 \le b_i \le M_2$ \cite{polyaszego1925}:
\begin{equation}
\frac{\left(\sum_{i=1}^n a_i^2\right) \left(\sum_{i=1}^n b_i^2\right)}{\left(\sum_{i=1}^n a_i b_i\right)^2} \le \frac{(M_1 M_2 + m_1 m_2)^2}{4 M_1 M_2 m_1 m_2}.
\end{equation}
Setting $b_i = 1$, $m_1 = 1$, and $M_1 = p$ recovers the equal-weight bound $I \le \frac{(p+1)^2}{4p}$.
\item \textbf{The Weighted Kantorovich Inequality (1948):} For arbitrary positive weights $w_i > 0$ and bounded positive ratios $0 < \mu_{\min} \le \mu_i \le \mu_{\max}$ with condition ratio $p = \mu_{\max}/\mu_{\min} \ge 1$ \cite{kantorovich1948}:
\begin{equation}
1 \le \frac{\left(\sum_{i=1}^n w_i \mu_i^2\right)\left(\sum_{i=1}^n w_i\right)}{\left(\sum_{i=1}^n w_i \mu_i\right)^2} \le \frac{(p+1)^2}{4p} = 1 + \frac{(p-1)^2}{4p}. \label{eq:kantorovich_weighted}
\end{equation}
Because the continuous relaxation parameter $\theta = \left(\sum_{i \in S_p} w_i\right) / \left(\sum_{i=1}^n w_i\right) \in [0, 1]$ can take arbitrary values in $[0,1]$ as the branch weights $w_i$ vary, the universal bound $\eta_{\mathrm{univ}}(p) = \frac{(p-1)^2}{4p}$ is sharp over all positively weighted quadratic forms. In Section~\ref{subsec:distortion_bridge}, we show that \eqref{eq:kantorovich_weighted} is the exact single-step squared-norm distortion of the linearization operator.
\end{enumerate}
\end{remark}

\section{The Method of Linearization}

The quadratic norm formulation in this section is defined on each flow-direction sign orthant and updated dynamically across sign transitions. More precisely, at iterate $\mathbf{x}^{(k)}$, let $\boldsymbol{\sigma}^{(k)}=\operatorname{sgn}(\mathbf{x}^{(k)})\in\{-1,1\}^n$ (with $\sigma_j^{(k)}=+1$ if $x_j^{(k)}=0$) partition each cycle $i$ into counter-directional branch sets $L_i^{(k)}$ and $R_i^{(k)}$. Once the iterate enters the open sign orthant of the equilibrium $\mathbf{x}^*$, the partition remains fixed and the local analysis applies directly on that open cone.

\subsection{Quadratic Form Representation}

Consider a nonzero quadratic form $\sum_{i=1}^s f_i^2$ with $f_1\ne0$. Defining relative coefficients $m_i = f_i / f_1$ with $m_1 = 1$:
\begin{equation}
\sum_{i=1}^s f_i^2 = \left( \frac{\sum_{i=1}^s m_i f_i}{\sqrt{1 + \sum_{i=2}^s m_i^2}} \right)^2. \label{eq:quadidentity}
\end{equation}
\begin{proof}
Direct algebraic expansion shows:
\begin{equation*}
\frac{\left(\sum_{i=1}^s m_i f_i\right)^2}{1 + \sum_{i=2}^s m_i^2} = \frac{\left(f_1 + \sum_{i=2}^s \frac{f_i^2}{f_1}\right)^2}{1 + \sum_{i=2}^s \frac{f_i^2}{f_1^2}} = \frac{\frac{1}{f_1^2}\left(\sum_{i=1}^s f_i^2\right)^2}{\frac{1}{f_1^2}\sum_{i=1}^s f_i^2} = \sum_{i=1}^s f_i^2.
\end{equation*}
Taking square roots yields the exact identity.
\end{proof}

\subsection{Linearization of a System of Equations}

Let a system of $c$ quadratic cycle equations in $n$ branch variables be partitioned by sign into counter-directional sets:
\begin{equation}
\sum_{k \in L_i} b_{ik} x_k^2 - \sum_{k \in R_i} |b_{ik}| x_k^2 = 0, \quad i = 1, 2, \dots, c, \label{eq:general_quad_sys}
\end{equation}
where $b_{ik} > 0$ for $k \in L_i$ and $b_{ik} < 0$ for $k \in R_i$. Equating both positive quadratic forms:
\begin{equation}
\sum_{k \in L_i} b_{ik} x_k^2 = \sum_{k \in R_i} |b_{ik}| x_k^2.
\end{equation}
Taking the square root of both sides gives the exact norm equality:
\begin{equation}
\sqrt{\sum_{k \in L_i} b_{ik} x_k^2} = \sqrt{\sum_{k \in R_i} |b_{ik}| x_k^2}.
\end{equation}

Applying identity \eqref{eq:quadidentity} to both positive forms around an operating point $\mathbf{x}^{(k)}$, we obtain the following linear relation for the next iterate $\mathbf{x}^{(k+1)}$:
\begin{equation}
\frac{\sum_{k \in L_i} b_{ik} x_k^{(k)} x_k^{(k+1)}}{\sqrt{\sum_{k \in L_i} b_{ik} (x_k^{(k)})^2}} = \frac{\sum_{k \in R_i} |b_{ik}| x_k^{(k)} x_k^{(k+1)}}{\sqrt{\sum_{k \in R_i} |b_{ik}| (x_k^{(k)})^2}}, \quad i = 1, \dots, c. \label{eq:linearized_sys}
\end{equation}

Coupling \eqref{eq:linearized_sys} with the linear junction continuity equations forms a complete coupled linear system solved simultaneously at each iteration. For a connected network represented by the square residual in Section~6, these are the reduced equations $\mathbf{A}_r\mathbf{x}^{(k+1)}=\mathbf{q}_r$.

\subsection{Single-Step Linearization Distortion via the Weighted Kantorovich Bound}
\label{subsec:distortion_bridge}

We now establish the exact mathematical link between the relative variance bound $\eta_{\mathrm{univ}}(p) = \frac{(p-1)^2}{4p}$ derived in Section~2 and the linearization step \eqref{eq:linearized_sys}.

In equation \eqref{eq:linearized_sys}, each exact directional norm $u_i(\mathbf{x}) = \sqrt{\sum_{j \in L_i} b_{ij} x_j^2}$ at a target flow vector $\mathbf{x}$ is approximated around the current operating point $\mathbf{x}^{(k)}$ by the linear functional
\begin{equation}
\widetilde{u}_i(\mathbf{x}; \mathbf{x}^{(k)}) \equiv \frac{\sum_{j \in L_i} b_{ij} x_j^{(k)} x_j}{\sqrt{\sum_{j \in L_i} b_{ij} (x_j^{(k)})^2}}. \label{eq:linear_surrogate}
\end{equation}

\begin{theorem}[Single-Step Squared-Norm Distortion Bound]
\label{thm:distortion_bridge}
Let $\mathbf{x}^{(k)}$ and $\mathbf{x}$ have strictly positive components on the directional index set $L_i$, and define the multiplicative coordinate ratios and directional condition ratio:
\begin{equation}
\mu_j = \frac{x_j}{x_j^{(k)}}, \quad \mu_{\min} = \min_{j \in L_i} \mu_j, \quad \mu_{\max} = \max_{j \in L_i} \mu_j, \quad p_i = \frac{\mu_{\max}}{\mu_{\min}} \ge 1.
\end{equation}
Then the linear surrogate $\widetilde{u}_i(\mathbf{x}; \mathbf{x}^{(k)})$ satisfies the two-sided sandwich inequality
\begin{equation}
1 \le \frac{u_i(\mathbf{x})^2}{\widetilde{u}_i(\mathbf{x}; \mathbf{x}^{(k)})^2} \le \frac{(p_i + 1)^2}{4p_i}, \label{eq:sandwich_bound}
\end{equation}
and the relative squared-norm linearization distortion is bounded sharply by $\eta_{\mathrm{univ}}(p_i)$:
\begin{equation}
0 \le \eta_{\mathrm{lin}}(\mathbf{x}; \mathbf{x}^{(k)}) \equiv \frac{u_i(\mathbf{x})^2 - \widetilde{u}_i(\mathbf{x}; \mathbf{x}^{(k)})^2}{\widetilde{u}_i(\mathbf{x}; \mathbf{x}^{(k)})^2} \le \frac{(p_i - 1)^2}{4p_i} = \eta_{\mathrm{univ}}(p_i), \label{eq:distortion_eta}
\end{equation}
with $\eta_{\mathrm{lin}}(\mathbf{x}; \mathbf{x}^{(k)}) = 0$ if and only if $\mathbf{x}$ is radially collinear with $\mathbf{x}^{(k)}$ on $L_i$ ($p_i = 1$).
\end{theorem}

\begin{proof}
Define the strictly positive branch weights $w_j = b_{ij} (x_j^{(k)})^2 > 0$ for $j \in L_i$. Substituting $x_j = \mu_j x_j^{(k)}$ into $u_i(\mathbf{x})^2$ and $\widetilde{u}_i(\mathbf{x}; \mathbf{x}^{(k)})^2$ gives
\begin{equation}
u_i(\mathbf{x})^2 = \sum_{j \in L_i} w_j \mu_j^2, \qquad \widetilde{u}_i(\mathbf{x}; \mathbf{x}^{(k)})^2 = \frac{\left(\sum_{j \in L_i} w_j \mu_j\right)^2}{\sum_{j \in L_i} w_j}.
\end{equation}
Taking their ratio yields the weighted relative-variance functional:
\begin{equation}
\frac{u_i(\mathbf{x})^2}{\widetilde{u}_i(\mathbf{x}; \mathbf{x}^{(k)})^2} = \frac{\left(\sum_{j \in L_i} w_j \mu_j^2\right)\left(\sum_{j \in L_i} w_j\right)}{\left(\sum_{j \in L_i} w_j \mu_j\right)^2} = 1 + \eta_{\mathrm{lin}}(\mathbf{x}; \mathbf{x}^{(k)}).
\end{equation}
The lower bound $1 \le u_i(\mathbf{x})^2 / \widetilde{u}_i(\mathbf{x}; \mathbf{x}^{(k)})^2$ follows immediately from the Cauchy--Schwarz inequality, with equality if and only if all $\mu_j$ are identical ($p_i = 1$). The upper bound follows directly from the weighted Kantorovich inequality \eqref{eq:kantorovich_weighted} in Remark~\ref{rem:kantorovich}. Subtracting $1$ across the inequality yields \eqref{eq:distortion_eta}.
\end{proof}

\begin{corollary}[Quadratic Contraction of the Multiplicative Distortion Envelope]
\label{cor:scalar_rec}
If the worst-case multiplicative squared-norm distortion factor $I(p_s) = 1 + \eta_{\mathrm{univ}}(p_s)$ at stage $s$ defines the subsequent condition ratio via the scalar envelope recurrence
\begin{equation}
p_{s+1} = I(p_s) = \frac{(p_s + 1)^2}{4p_s}, \quad p_1 > 1, \label{eq:scalar_rec}
\end{equation}
then $p_{s+1} - 1 = \frac{(p_s - 1)^2}{4p_s} > 0$ and $p_s - p_{s+1} = \frac{(p_s - 1)(3p_s + 1)}{4p_s} > 0$. Consequently, $\{p_s\}$ decreases monotonically to $1$ with a quadratic asymptotic convergence rate:
\begin{equation}
\lim_{s \to \infty} \frac{p_{s+1} - 1}{(p_s - 1)^2} = \frac{1}{4}.
\end{equation}
\end{corollary}

\begin{remark}[Scope of the Scalar Envelope Recurrence]
Corollary~\ref{cor:scalar_rec} characterizes the idealized scalar contraction of the worst-case Kantorovich distortion factor $I(p)$ on an isolated directional norm. In a coupled network, the iterate $\mathbf{x}^{(k+1)}$ is determined jointly across all $c$ cycles and $m-1$ continuity constraints, and its rigorous multidimensional convergence rate is governed by the Newton--Kantorovich analysis in Section~4.
\end{remark}

\section{Rigorous Convergence Proofs and Nonlinear Preconditioning}

In this section, we establish the algebraic equivalence between the linearized system \eqref{eq:linearized_sys} and a multidimensional Newton--Raphson iteration applied to norm differences. We then prove that the coupled Jacobian is nonsingular at every regular equilibrium, establish local quadratic convergence, and show how the degree-one homogeneous norm formulation acts as a nonlinear degree-reducing preconditioner relative to standard Newton on raw polynomial loop equations.

\subsection{Equivalence to Newton--Raphson on Norm Differences}

\begin{theorem}[Equivalence Theorem]
\label{thm:newton_equivalence}
The linearized equation \eqref{eq:linearized_sys} is algebraically identical to the classical Newton--Raphson method applied to the difference of Euclidean norms:
\begin{equation}
F_i(\mathbf{x}) \equiv u_i(\mathbf{x}) - v_i(\mathbf{x}) = 0, \label{eq:norm_diff}
\end{equation}
where $u_i(\mathbf{x}) = \sqrt{\sum_{k \in L_i} b_{ik} x_k^2}$ and $v_i(\mathbf{x}) = \sqrt{\sum_{k \in R_i} |b_{ik}| x_k^2}$.
When linear continuity constraints are appended, their Newton rows coincide with the constraints themselves, so the full coupled square system is likewise the Newton system for the combined residual.
\end{theorem}

\begin{proof}
Let $\mathbf{x}$ lie in a fixed-direction open region in which the two directional norms in each equation are nonzero. The directional Euclidean norm $u_i(\mathbf{x})$ is positively homogeneous of degree 1:
\begin{equation*}
u_i(\lambda \mathbf{x}) = \lambda u_i(\mathbf{x}) \quad \forall \lambda > 0.
\end{equation*}
By Euler's Homogeneous Function Theorem:
\begin{equation}
\nabla u_i(\mathbf{x})^T \mathbf{x} = u_i(\mathbf{x}), \qquad \nabla v_i(\mathbf{x})^T \mathbf{x} = v_i(\mathbf{x}). \label{eq:euler}
\end{equation}
Subtracting these relations gives:
\begin{equation}
\nabla F_i(\mathbf{x})^T \mathbf{x} = u_i(\mathbf{x}) - v_i(\mathbf{x}) = F_i(\mathbf{x}). \label{eq:euler_diff}
\end{equation}

Now consider the classical multidimensional Newton--Raphson step for solving $F_i(\mathbf{x}) = 0$ linearized around iterate $\mathbf{x}^{(k)}$:
\begin{equation}
\nabla F_i(\mathbf{x}^{(k)})^T \left(\mathbf{x}^{(k+1)} - \mathbf{x}^{(k)}\right) = -F_i(\mathbf{x}^{(k)}).
\end{equation}
Rearranging to isolate $\mathbf{x}^{(k+1)}$ on the left-hand side:
\begin{equation}
\nabla F_i(\mathbf{x}^{(k)})^T \mathbf{x}^{(k+1)} = \nabla F_i(\mathbf{x}^{(k)})^T \mathbf{x}^{(k)} - F_i(\mathbf{x}^{(k)}). \label{eq:newton_rearranged}
\end{equation}
Substituting Euler's identity \eqref{eq:euler_diff} into the right-hand side of \eqref{eq:newton_rearranged}:
\begin{equation}
\nabla F_i(\mathbf{x}^{(k)})^T \mathbf{x}^{(k+1)} = F_i(\mathbf{x}^{(k)}) - F_i(\mathbf{x}^{(k)}) \equiv 0.
\end{equation}
Remarkably, \textbf{the affine constant term vanishes identically}. 

Expanding $\nabla F_i = \nabla u_i - \nabla v_i$:
\begin{equation}
\nabla u_i(\mathbf{x}^{(k)})^T \mathbf{x}^{(k+1)} = \nabla v_i(\mathbf{x}^{(k)})^T \mathbf{x}^{(k+1)}. \label{eq:grad_balance}
\end{equation}
Computing the partial derivatives explicitly:
\begin{equation*}
\frac{\partial u_i}{\partial x_j} = \frac{b_{ij} x_j}{\sqrt{\sum_{k \in L_i} b_{ik} x_k^2}}, \qquad \frac{\partial v_i}{\partial x_j} = \frac{|b_{ij}| x_j}{\sqrt{\sum_{k \in R_i} |b_{ik}| x_k^2}}.
\end{equation*}
Substituting these into \eqref{eq:grad_balance} yields:
\begin{equation}
\frac{\sum_{k \in L_i} b_{ik} x_k^{(k)} x_k^{(k+1)}}{\sqrt{\sum_{k \in L_i} b_{ik} (x_k^{(k)})^2}} = \frac{\sum_{k \in R_i} |b_{ik}| x_k^{(k)} x_k^{(k+1)}}{\sqrt{\sum_{k \in R_i} |b_{ik}| (x_k^{(k)})^2}},
\end{equation}
which is exactly equation \eqref{eq:linearized_sys}.
\end{proof}

\subsection{Equilibrium Row-Equilibration Identity and Quadratic Convergence}

\begin{proposition}[Equilibrium Row-Equilibration and Nonsingularity of $\mathbf{J}_{\mathbf{F}}(\mathbf{x}^*)$]
\label{prop:nonsingularity}
Let $\mathbf{x}^*$ be an equilibrium of a connected network with reduced incidence matrix $\mathbf{A}_r \in \mathbb{R}^{(m-1)\times n}$, fundamental cycle matrix $\mathbf{C} \in \mathbb{R}^{c \times n}$, resistance vector $\mathbf{s} \in \mathbb{R}_{++}^n$, and exponent $\alpha > 1$, such that each cycle has nonzero directional norm $w_i \equiv u_i(\mathbf{x}^*) = v_i(\mathbf{x}^*) > 0$ ($i = 1, \dots, c$). Define the diagonal matrices
\begin{equation}
\mathbf{W}(\mathbf{x}^*) = \operatorname{diag}(w_1, \dots, w_c) \succ \mathbf{0}, \qquad \mathbf{D}(\mathbf{x}^*) = \operatorname{diag}\!\left(s_1 |x_1^*|^{\alpha-1}, \dots, s_n |x_n^*|^{\alpha-1}\right) \succeq \mathbf{0}.
\end{equation}
Then the loop Jacobian block $\mathbf{J}_{\text{loop}}(\mathbf{x}^*)$ of the norm-difference residual satisfies the exact row-equilibration identity
\begin{equation}
\mathbf{J}_{\text{loop}}(\mathbf{x}^*) = \mathbf{W}(\mathbf{x}^*)^{-(\alpha-1)} \mathbf{C} \mathbf{D}(\mathbf{x}^*). \label{eq:row_equilibration}
\end{equation}
Furthermore, the full square Jacobian $\mathbf{J}_{\mathbf{F}}(\mathbf{x}^*) = \begin{bmatrix} \mathbf{A}_r \\ \mathbf{J}_{\text{loop}}(\mathbf{x}^*) \end{bmatrix} \in \mathbb{R}^{n \times n}$ is nonsingular if and only if $\mathbf{C} \mathbf{D}(\mathbf{x}^*) \mathbf{C}^T \succ \mathbf{0}$ (which holds whenever the set of branches with nonzero equilibrium flow contains a spanning tree of $\mathcal{G}$, and in particular whenever $x_k^* \ne 0$ for all $k$).
\end{proposition}

\begin{proof}
At $\mathbf{x}^*$, $u_i(\mathbf{x}^*) = v_i(\mathbf{x}^*) = w_i > 0$. For any branch $k$ in cycle $i$, $C_{ik} \in \{-1, 1\}$ and $\sigma_k^* = \operatorname{sgn}(x_k^*)$. If $k \in L_i$, then $C_{ik} \sigma_k^* = +1 \implies C_{ik} = \sigma_k^*$, so
\begin{equation*}
[\mathbf{J}_{\text{loop}}(\mathbf{x}^*)]_{ik} = \frac{s_k |x_k^*|^{\alpha-1} \sigma_k^*}{u_i(\mathbf{x}^*)^{\alpha-1}} = \frac{C_{ik} s_k |x_k^*|^{\alpha-1}}{w_i^{\alpha-1}}.
\end{equation*}
If $k \in R_i$, then $C_{ik} \sigma_k^* = -1 \implies C_{ik} = -\sigma_k^*$, and therefore
\begin{equation*}
[\mathbf{J}_{\text{loop}}(\mathbf{x}^*)]_{ik} = -\frac{s_k |x_k^*|^{\alpha-1} \sigma_k^*}{v_i(\mathbf{x}^*)^{\alpha-1}} = \frac{C_{ik} s_k |x_k^*|^{\alpha-1}}{w_i^{\alpha-1}},
\end{equation*}
which proves \eqref{eq:row_equilibration}. Now suppose $\mathbf{J}_{\mathbf{F}}(\mathbf{x}^*) \mathbf{v} = \mathbf{0}$ for some $\mathbf{v} \in \mathbb{R}^n$. Then $\mathbf{A}_r \mathbf{v} = \mathbf{0}$, so $\mathbf{v} \in \ker(\mathbf{A}_r) = \operatorname{range}(\mathbf{C}^T)$, which implies $\mathbf{v} = \mathbf{C}^T \boldsymbol{\beta}$ for a unique $\boldsymbol{\beta} \in \mathbb{R}^c$. Substituting into the bottom block gives
\begin{equation*}
\mathbf{W}(\mathbf{x}^*)^{-(\alpha-1)} \mathbf{C} \mathbf{D}(\mathbf{x}^*) \mathbf{C}^T \boldsymbol{\beta} = \mathbf{0} \iff \mathbf{C} \mathbf{D}(\mathbf{x}^*) \mathbf{C}^T \boldsymbol{\beta} = \mathbf{0}.
\end{equation*}
Thus $\mathbf{J}_{\mathbf{F}}(\mathbf{x}^*)$ is nonsingular if and only if $\mathbf{C} \mathbf{D}(\mathbf{x}^*) \mathbf{C}^T \succ \mathbf{0}$.
\end{proof}

\begin{theorem}[Local Quadratic Convergence]
\label{thm:quad_conv}
Let $\mathbf{x}^*$ be an equilibrium in a fixed-direction open region with $x_k^* \ne 0$ for all $k$ and $\alpha \ge 2$ (or $\alpha > 1$ with strictly nonzero branch flows). Then $\mathbf{F}$ is $C^2$ in a neighborhood of $\mathbf{x}^*$, $\mathbf{J}_{\mathbf{F}}(\mathbf{x}^*)$ is nonsingular by Proposition~\ref{prop:nonsingularity}, and there exists an open neighborhood $B(\mathbf{x}^*,r)$ such that for any $\mathbf{x}^{(0)}\in B(\mathbf{x}^*,r)$, the iterates generated by \eqref{eq:linearized_sys} converge to $\mathbf{x}^*$ quadratically:
\begin{equation}
\|\mathbf{x}^{(k+1)} - \mathbf{x}^*\| \le C \|\mathbf{x}^{(k)} - \mathbf{x}^*\|^2.
\end{equation}
\end{theorem}

\begin{proof}
By Theorem~\ref{thm:newton_equivalence}, the iteration is the Newton map $\mathbf{x}^{(k+1)}=\mathbf{x}^{(k)}-\mathbf{J}_{\mathbf{F}}(\mathbf{x}^{(k)})^{-1}\mathbf{F}(\mathbf{x}^{(k)})$. By Proposition~\ref{prop:nonsingularity}, $\mathbf{J}_{\mathbf{F}}(\mathbf{x}^*)$ is nonsingular and $\mathbf{J}_{\mathbf{F}}$ is locally Lipschitz on $B(\mathbf{x}^*,r)$. The classical Newton--Kantorovich theorem \cite{kantorovich1948,vasiliev1980} therefore yields local quadratic convergence.
\end{proof}

\subsection{Nonlinear Degree-Reducing Preconditioning (\texorpdfstring{$u_i^d - v_i^d \mapsto u_i - v_i$}{u\_i\^d - v\_i\^d to u\_i - v\_i})}
\label{subsec:preconditioning}

Although Theorem~\ref{thm:newton_equivalence} establishes that the scheme is a Newton iteration on $F_i(\mathbf{x}) = u_i(\mathbf{x}) - v_i(\mathbf{x}) = 0$, standard Newton implementations for hydraulic and electrical networks \cite{todini1987} apply Newton's method directly to the raw degree-$d$ head-loss residual ($d=2$ for quadratic Darcy--Weisbach systems; $d=\alpha>1$ for general power laws):
\begin{equation}
\widetilde{F}_i(\mathbf{x}) \equiv u_i(\mathbf{x})^d - v_i(\mathbf{x})^d = \sum_{k \in L_i} b_{ik} |x_k|^d - \sum_{k \in R_i} |b_{ik}| |x_k|^d = 0, \label{eq:raw_poly_res}
\end{equation}
where $u_i(\mathbf{x}) = \left(\sum_{k \in L_i} b_{ik} |x_k|^d\right)^{1/d}$ and $v_i(\mathbf{x}) = \left(\sum_{k \in R_i} |b_{ik}| |x_k|^d\right)^{1/d}$.

On the open domain where $u_i(\mathbf{x}) > 0$ and $v_i(\mathbf{x}) > 0$, both residuals share the exact same zero manifold: $\widetilde{F}_i(\mathbf{x}) = 0 \iff F_i(\mathbf{x}) = 0$. However, replacing $\widetilde{F}_i(\mathbf{x})$ with $F_i(\mathbf{x})$ reduces the homogeneity degree of the residual from $d$ to $1$. This transformation acts as a \textbf{nonlinear degree-reducing preconditioner} with four fundamental numerical advantages over standard Newton on $\widetilde{F}_i(\mathbf{x}) = 0$:

\begin{proposition}[Structural Properties of Degree-One Norm Preconditioning]
\label{prop:preconditioning}
Let $d > 1$ (with $d=2$ for quadratic Darcy--Weisbach networks and $d=\alpha>1$ for general power-law systems) and let $\mathbf{x}$ belong to an open fixed-direction region where $u_i(\mathbf{x}), v_i(\mathbf{x}) > 0$ and $u_i, v_i$ are $C^2$. Comparing the degree-one residual $F_i(\mathbf{x}) = u_i(\mathbf{x}) - v_i(\mathbf{x})$ against the raw degree-$d$ residual $\widetilde{F}_i(\mathbf{x}) = u_i(\mathbf{x})^d - v_i(\mathbf{x})^d$:
\begin{enumerate}
\item \textbf{Elimination of Radial Curvature:} Differentiating Euler's identity $\nabla u_i(\mathbf{x})^T \mathbf{x} = u_i(\mathbf{x})$ with respect to $\mathbf{x}$ yields $\nabla^2 u_i(\mathbf{x})\mathbf{x} + \nabla u_i(\mathbf{x}) = \nabla u_i(\mathbf{x})$, and therefore
\begin{equation}
\nabla^2 u_i(\mathbf{x})\mathbf{x} = \mathbf{0}, \qquad \nabla^2 v_i(\mathbf{x})\mathbf{x} = \mathbf{0}, \qquad \nabla^2 F_i(\mathbf{x})\mathbf{x} = \mathbf{0}. \label{eq:radial_hessian_zero}
\end{equation}
The Hessian of $F_i$ vanishes identically along radial rays. By contrast, because $\widetilde{F}_i$ is positively homogeneous of degree $d$, its radial curvature is nonzero everywhere outside the solution manifold:
\begin{equation}
\nabla^2 \widetilde{F}_i(\mathbf{x})\mathbf{x} = (d - 1)\nabla \widetilde{F}_i(\mathbf{x}) \ne \mathbf{0}. \label{eq:radial_hessian_poly}
\end{equation}
\item \textbf{Cancellation of the Affine Residual Shift:} For the degree-$d$ residual $\widetilde{F}_i$, Euler's identity gives $\nabla \widetilde{F}_i(\mathbf{x}^{(k)})^T \mathbf{x}^{(k)} = d \widetilde{F}_i(\mathbf{x}^{(k)})$, so the rearranged Newton equation retains a non-vanishing affine right-hand side:
\begin{equation}
\nabla \widetilde{F}_i(\mathbf{x}^{(k)})^T \mathbf{x}^{(k+1)} = (d - 1)\widetilde{F}_i(\mathbf{x}^{(k)}). \label{eq:poly_newton_rhs}
\end{equation}
When $\mathbf{x}^{(k)}$ is far from equilibrium, $(d-1)\widetilde{F}_i(\mathbf{x}^{(k)})$ introduces large affine shifts that delay convergence under fixed step budgets or require line-search damping. Reducing the homogeneity degree to $d=1$ makes the factor $(d-1)$ vanish identically, yielding the unshifted linear conservation law $\nabla F_i(\mathbf{x}^{(k)})^T \mathbf{x}^{(k+1)} = 0$.
\item \textbf{Degree-Zero Jacobian Scale Invariance:} The gradient $\nabla F_i(\mathbf{x})$ is positively homogeneous of degree $0$:
\begin{equation}
\nabla F_i(\lambda \mathbf{x}) = \nabla F_i(\mathbf{x}) \quad \forall \lambda > 0,
\end{equation}
and its $\ell_{d/(d-1)}$-weighted norm is bounded independently of flow magnitude. Conversely, $\nabla \widetilde{F}_i(\lambda \mathbf{x}) = \lambda^{d-1}\nabla \widetilde{F}_i(\mathbf{x})$ scales polynomially with flow magnitude.
\item \textbf{Automatic Equilibrium Row-Equilibration:} By Proposition~\ref{prop:nonsingularity}, at equilibrium $\mathbf{J}_{\text{loop}}(\mathbf{x}^*) = \frac{1}{d}\mathbf{W}(\mathbf{x}^*)^{-(d-1)}\mathbf{J}_{\widetilde{\mathbf{F}},\text{loop}}(\mathbf{x}^*)$, which automatically divides each loop row $i$ by its characteristic head-loss scale $w_i^{d-1}$ across multi-scale circuits.
\end{enumerate}
\end{proposition}

\begin{remark}[Effect on Newton Basin of Attraction]
Because $\nabla F_i(\lambda \mathbf{x}^*) = \nabla F_i(\mathbf{x}^*)$ and $\nabla^2 F_i(\mathbf{x})\mathbf{x} = \mathbf{0}$, any radial scaling error in the initialization $\mathbf{x}^{(0)} = \lambda \mathbf{x}^*$ produces zero distortion in the linearized loop equations ($\eta_{\mathrm{lin}} = 0$ in Theorem~\ref{thm:distortion_bridge}), allowing the coupled continuity-loop system to recover the exact scale in a single linear solve. Under standard Newton on $\widetilde{F}_i(\mathbf{x}) = 0$, a pure scaling error contracts only linearly by the factor $\frac{d-1}{d}$ per step (halving the error each step for quadratic $d=2$). This directly explains why the norm-linearized solver reaches tolerance within a fixed $K=4$ budget where undamped polynomial Newton fails on $100\%$ of cold-start initializations (Section~\ref{subsec:experiments}).
\end{remark}

\section{Generalizations to Non-Quadratic Power Laws and Modern Optimization}

\subsection{Extension to General Power Laws via Homogeneity}

In networked physical systems, energy dissipation can be governed by non-integer power laws $h_i=s_i|x_i|^{\alpha-1}x_i$ (for example, the Hazen--Williams exponent is $\alpha=1.852$ in hydraulic pipe networks \cite{todini1987}). Under dynamic sign partitioning $\boldsymbol{\sigma}^{(k)}=\operatorname{sgn}(\mathbf{x}^{(k)})$ of a cycle balance, the magnitude equation can be written as
\begin{equation}
\sum_{k \in L_i} b_{ik} |x_k|^\alpha = \sum_{k \in R_i} |b_{ik}| |x_k|^\alpha.
\end{equation}
Let $p=\alpha>1$, and define the weighted $\ell_p$ norm:
\begin{equation}
u_i(\mathbf{x}) = \left(\sum_{k \in L_i} b_{ik} |x_k|^p\right)^{1/p}.
\end{equation}
On an open region where $u_i$ and $v_i$ are nonzero, this norm is differentiable for $p>1$ and positively homogeneous of degree 1, so Euler's theorem gives $\nabla u_i(\mathbf{x})^T\mathbf{x}=u_i(\mathbf{x})$. Computing the gradient explicitly (with the numerator interpreted as zero at $x_j=0$):
\begin{equation}
\frac{\partial u_i}{\partial x_j} = \frac{b_{ij} |x_j|^{p-1} \operatorname{sgn}(x_j)}{\left(\sum_{k \in L_i} b_{ik} |x_k|^p\right)^{\frac{p-1}{p}}}.
\end{equation}
Equating the first-order projections $\nabla u_i(\mathbf{x}^{(k)})^T \mathbf{x}^{(k+1)} = \nabla v_i(\mathbf{x}^{(k)})^T \mathbf{x}^{(k+1)}$ yields the \textbf{generalized linear theory}:
\begin{equation}
\frac{\sum_{k \in L_i} b_{ik} |x_k^{(k)}|^{p-1}\operatorname{sgn}(x_k^{(k)})x_k^{(k+1)}}{\left(\sum_{k \in L_i} b_{ik} |x_k^{(k)}|^p\right)^{\frac{p-1}{p}}} = \frac{\sum_{k \in R_i} |b_{ik}| |x_k^{(k)}|^{p-1}\operatorname{sgn}(x_k^{(k)})x_k^{(k+1)}}{\left(\sum_{k \in R_i} |b_{ik}| |x_k^{(k)}|^p\right)^{\frac{p-1}{p}}}.
\end{equation}
Thus the vanishing affine-term identity and the degree-reducing preconditioning properties of Proposition~\ref{prop:preconditioning} extend directly to weighted $\ell_p$ norm differences for $p>1$. The case $p=1$ is linear and solved in a single step, which provides a natural cold-start initialization $\mathbf{x}_{\text{lin}}^{(0)}$ for $p>1$.

\subsection{Difference-of-Convex (DC) Programming Context}

The equation $u_i(\mathbf{x})-v_i(\mathbf{x})=0$ is the zero set of a difference of two convex functions. While this representation is a root-finding problem rather than a scalar convex optimization problem, convexity of the directional norms $u_i$ and $v_i$ implies that their first-order Taylor approximations are global lower-bounding supporting hyperplanes: $u_i(\mathbf{x}) \ge \nabla u_i(\mathbf{x}^{(k)})^T \mathbf{x}$, which corresponds directly to the lower bound $1 \le u_i(\mathbf{x})^2 / \widetilde{u}_i(\mathbf{x}; \mathbf{x}^{(k)})^2$ in Theorem~\ref{thm:distortion_bridge}.

\section{Differentiable Physics Implicit Layers (Deep Equilibrium Models)}

We now deploy the theoretical results established above into modern deep learning architectures as a \textbf{Differentiable Physics Implicit Layer (Deep Equilibrium Model / DEQ)}.

\subsection{Mathematical Architecture and Global Uniqueness}

Let $\mathcal{G}=(\mathcal{V},\mathcal{E},\mathcal{L})$ be a connected directed physical network graph with $m$ junction nodes, $n$ branches (edges), and $c=n-m+1$ fundamental independent cycles. Let $\mathbf{A}_r\in\mathbb{R}^{(m-1)\times n}$ be a reduced incidence matrix obtained by deleting one redundant row of the full incidence matrix, and let $\mathbf{q}_r\in\mathbb{R}^{m-1}$ be the associated reduced demand vector.

\begin{proposition}[Global Existence and Uniqueness of the Network Equilibrium]
\label{prop:global_uniqueness}
For any resistance vector $\mathbf{s}\in\mathbb{R}_{++}^n$, exponent $\alpha > 1$, and reduced demand $\mathbf{q}_r\in\mathbb{R}^{m-1}$, the strictly convex, coercive dissipation potential minimization problem
\begin{equation}
\min_{\mathbf{x}\in\mathbb{R}^n}\;\Phi(\mathbf{x};\mathbf{s}) \equiv \sum_{k=1}^n \frac{s_k}{\alpha+1}|x_k|^{\alpha+1} \quad \text{subject to} \quad \mathbf{A}_r\mathbf{x}=\mathbf{q}_r \label{eq:variational_potential}
\end{equation}
admits a unique global minimizer $\mathbf{x}^*(\mathbf{s},\mathbf{q}_r)\in\mathbb{R}^n$, which is the unique solution of the coupled conservation system $\mathbf{A}_r\mathbf{x}^*=\mathbf{q}_r$ and $u_i(\mathbf{x}^*)-v_i(\mathbf{x}^*)=0$ ($i=1,\dots,c$).
\end{proposition}

\begin{proof}
Because $\mathbf{A}_r$ has full row rank $m-1$, the affine feasible set $\{\mathbf{x}\in\mathbb{R}^n : \mathbf{A}_r\mathbf{x}=\mathbf{q}_r\}$ is nonempty, closed, and convex. Since $s_k > 0$ and $\alpha > 1$, $\Phi(\cdot;\mathbf{s})$ is strictly convex, $C^1$, and coercive ($\Phi(\mathbf{x};\mathbf{s})\to\infty$ as $\|\mathbf{x}\|\to\infty$), guaranteeing a unique global minimizer $\mathbf{x}^*$. By first-order KKT optimality, $\mathbf{x}^*$ is the minimizer if and only if $\nabla_{\mathbf{x}}\Phi(\mathbf{x}^*;\mathbf{s}) = \mathbf{s}\odot|\mathbf{x}^*|^{\alpha-1}\odot\mathbf{x}^* \in \operatorname{range}(\mathbf{A}_r^T) = \ker(\mathbf{C})$, which is equivalent to $\mathbf{C}(\mathbf{s}\odot|\mathbf{x}^*|^{\alpha-1}\odot\mathbf{x}^*)=\mathbf{0}$ and hence $u_i(\mathbf{x}^*)-v_i(\mathbf{x}^*)=0$ for all $i=1,\dots,c$.
\end{proof}

\begin{definition}[Differentiable Network Flow Layer]
Given learnable physical parameters $\boldsymbol{\theta}=(\mathbf{s},\mathbf{q}_r)\in\mathbb{R}_{++}^n\times\mathbb{R}^{m-1}$, the \textbf{FlowDEQ Layer} maps $\boldsymbol{\theta}$ to the unique global equilibrium flow state $\mathbf{x}^*(\boldsymbol{\theta})\in\mathbb{R}^n$ of Proposition~\ref{prop:global_uniqueness}:
\begin{equation}
\mathbf{x}^*(\boldsymbol{\theta}) = \operatorname{FlowDEQ}(\boldsymbol{\theta}) \iff \mathbf{F}(\mathbf{x}^*, \boldsymbol{\theta}) = \mathbf{0},
\end{equation}
where the square physical operator $\mathbf{F}:\mathbb{R}^n\times\mathbb{R}^{n+m-1}\to\mathbb{R}^n$ combines reduced junction continuity and cycle norm balances partitioned by $\boldsymbol{\sigma}=\operatorname{sgn}(\mathbf{x})$:
\begin{equation}
\mathbf{F}(\mathbf{x}, \boldsymbol{\theta}) = \begin{bmatrix}
\mathbf{A}_r \mathbf{x} - \mathbf{q}_r \\[1ex]
u_1(\mathbf{x}) - v_1(\mathbf{x}) \\
\vdots \\
u_c(\mathbf{x}) - v_c(\mathbf{x})
\end{bmatrix} = \mathbf{0}. \label{eq:deq_operator}
\end{equation}
\end{definition}

\subsection{Forward Pass: Dynamic-Sign Fixed-Budget Implementation}

In standard Deep Equilibrium Models \cite{bai2019}, generic fixed-point solvers such as Broyden's method can require instance-dependent numbers of iterations. In a batched GPU implementation, a fixed computational budget avoids dynamic branching across batch elements.

In our framework, a forward pass with a fixed budget $K$ executes the norm-linearization scheme without requiring oracle knowledge of flow directions:
\begin{enumerate}
\item Obtain an oracle-free cold-start initialization $\mathbf{x}^{(0)}$ from the linear resistive solve ($\alpha=1$):
\begin{equation}
\begin{bmatrix} \mathbf{A}_r \\ \mathbf{C}\operatorname{diag}(\mathbf{s}_{\text{prior}}) \end{bmatrix} \mathbf{x}_{\text{lin}} = \begin{bmatrix} \mathbf{q}_r \\ \mathbf{0} \end{bmatrix}.
\end{equation}
\item For $k=0,\ldots,K-1$, update the branch sign vector $\boldsymbol{\sigma}^{(k)}=\operatorname{sgn}(\mathbf{x}^{(k)})$ and assemble the linearized coupled matrix $\mathbf{M}(\mathbf{x}^{(k)}) \in \mathbb{R}^{n \times n}$:
\begin{equation}
\mathbf{M}(\mathbf{x}^{(k)}) \mathbf{x}^{(k+1)} = \begin{bmatrix} \mathbf{q}_r \\ \mathbf{0} \end{bmatrix},
\end{equation}
where the top $m-1$ rows enforce exact reduced nodal continuity $\mathbf{A}_r \mathbf{x}^{(k+1)} = \mathbf{q}_r$ and the bottom $c$ rows correspond to the scale-invariant, unshifted gradient projections \eqref{eq:linearized_sys}.
\item Solve the nonsingular linear system using a batched direct LU solver.
\end{enumerate}

In the experiments below, $K=4$ is selected as a fixed operational budget, paired with a post-solve residual verification against the target tolerance.

\subsection{Backward Pass: Closed-Form Adjoint Differentiation via the IFT}

During neural network training, a downstream scalar objective $\mathcal{L}(\mathbf{x}^*(\boldsymbol{\theta}))$ is evaluated. To update the parameters via gradient descent, we require the Vector-Jacobian Product (VJP) $\nabla_{\boldsymbol{\theta}} \mathcal{L} = \frac{d \mathcal{L}}{d \boldsymbol{\theta}}$.

\begin{theorem}[Exact IFT Adjoint and Closed-Form Parameter Gradients]
\label{thm:ift_adjoint}
Let $\mathbf{x}^*$ be a regular equilibrium with nonzero branch flows, so that $\mathbf{J}_{\mathbf{F}}(\mathbf{x}^*,\boldsymbol{\theta})$ is nonsingular by Proposition~\ref{prop:nonsingularity}. Partition the adjoint state vector $\boldsymbol{\lambda}^* = \begin{bmatrix} \boldsymbol{\lambda}_{\text{node}}^*;\; \boldsymbol{\lambda}_{\text{loop}}^* \end{bmatrix} \in \mathbb{R}^{(m-1)+c}$ solving the transposed linear system
\begin{equation}
\mathbf{J}_{\mathbf{F}}(\mathbf{x}^*, \boldsymbol{\theta})^T \boldsymbol{\lambda}^* = \nabla_{\mathbf{x}} \mathcal{L}(\mathbf{x}^*). \label{eq:adjoint_sys}
\end{equation}
Then the exact gradients of $\mathcal{L}$ with respect to $\mathbf{q}_r$ and $\mathbf{s}$ are given in closed form by:
\begin{equation}
\nabla_{\mathbf{q}_r} \mathcal{L} = \boldsymbol{\lambda}_{\text{node}}^*, \qquad \frac{\partial \mathcal{L}}{\partial s_k} = -\frac{1}{\alpha}|x_k^*|^\alpha \sum_{i=1}^c \lambda_{\text{loop},i}^* \left(\frac{\mathbb{I}(k\in L_i)}{u_i(\mathbf{x}^*)^{\alpha-1}} - \frac{\mathbb{I}(k\in R_i)}{v_i(\mathbf{x}^*)^{\alpha-1}}\right), \label{eq:vjp_closed}
\end{equation}
which for quadratic Darcy--Weisbach systems ($\alpha=2$) simplifies to
\begin{equation*}
\frac{\partial \mathcal{L}}{\partial s_k} = -\frac{1}{2}(x_k^*)^2 \sum_{i=1}^c \lambda_{\text{loop},i}^* \left(\frac{\mathbb{I}(k\in L_i)}{u_i(\mathbf{x}^*)} - \frac{\mathbb{I}(k\in R_i)}{v_i(\mathbf{x}^*)}\right).
\end{equation*}
\end{theorem}

\begin{proof}
Differentiating the equilibrium identity $\mathbf{F}(\mathbf{x}^*(\boldsymbol{\theta}), \boldsymbol{\theta}) \equiv \mathbf{0}$ with respect to $\boldsymbol{\theta}$ via the Implicit Function Theorem yields $\frac{\partial \mathbf{x}^*}{\partial \boldsymbol{\theta}} = - \mathbf{J}_{\mathbf{F}}(\mathbf{x}^*, \boldsymbol{\theta})^{-1} \frac{\partial \mathbf{F}(\mathbf{x}^*, \boldsymbol{\theta})}{\partial \boldsymbol{\theta}}$, and therefore
\begin{equation}
\nabla_{\boldsymbol{\theta}} \mathcal{L} = - \left( \frac{\partial \mathbf{F}(\mathbf{x}^*, \boldsymbol{\theta})}{\partial \boldsymbol{\theta}} \right)^T \boldsymbol{\lambda}^*. \label{eq:vjp}
\end{equation}
Since $\mathbf{F}_{\text{node}}(\mathbf{x},\boldsymbol{\theta}) = \mathbf{A}_r\mathbf{x} - \mathbf{q}_r$, we have $\frac{\partial \mathbf{F}_{\text{node}}}{\partial \mathbf{q}_r} = -\mathbf{I}_{m-1}$ and $\frac{\partial \mathbf{F}_{\text{loop}}}{\partial \mathbf{q}_r} = \mathbf{0}$, giving $\nabla_{\mathbf{q}_r}\mathcal{L} = \boldsymbol{\lambda}_{\text{node}}^*$. Similarly, differentiating $F_{\text{loop},i} = u_i(\mathbf{x}^*) - v_i(\mathbf{x}^*)$ with respect to $s_k$ yields
\begin{equation*}
\frac{\partial F_{\text{loop},i}}{\partial s_k} = \frac{1}{\alpha}|x_k^*|^\alpha \left(\frac{\mathbb{I}(k\in L_i)}{u_i(\mathbf{x}^*)^{\alpha-1}} - \frac{\mathbb{I}(k\in R_i)}{v_i(\mathbf{x}^*)^{\alpha-1}}\right),
\end{equation*}
which upon substitution into \eqref{eq:vjp} gives \eqref{eq:vjp_closed}.
\end{proof}

\begin{remark}[Approximate Solves, Memory, and Factorizations]
This adjoint formulation has the following implementation consequences:
\begin{enumerate}
\item \textbf{Memory with respect to solver steps:} An IFT implementation retains only the terminal state and local derivative data needed at the evaluation point, rather than unrolling a $K$-step computational graph, yielding $O(1)$ activation memory with respect to $K$.
\item \textbf{Approximate IFT gradient after early stopping:} If the forward pass stops at $\widetilde{\mathbf{x}}$ with a nonzero residual, then solving
\begin{equation}
\mathbf{J}_{\mathbf{F}}(\widetilde{\mathbf{x}},\boldsymbol{\theta})^T\widetilde{\boldsymbol{\lambda}}=\nabla_{\mathbf{x}}\mathcal{L}(\widetilde{\mathbf{x}}),
\end{equation}
and using $-\mathbf{F}_{\boldsymbol{\theta}}(\widetilde{\mathbf{x}},\boldsymbol{\theta})^T\widetilde{\boldsymbol{\lambda}}$ yields an approximate equilibrium gradient whose approximation error scales linearly with the terminal residual $\|\mathbf{F}(\widetilde{\mathbf{x}}, \boldsymbol{\theta})\|$.
\item \textbf{Factorizations:} The final forward Newton step from $\mathbf{x}^{(K-1)}$ factors $\mathbf{J}_{\mathbf{F}}(\mathbf{x}^{(K-1)})$, whereas the exact adjoint at $\mathbf{x}^{(K)}$ evaluates $\mathbf{J}_{\mathbf{F}}(\mathbf{x}^{(K)})^T$. When $\mathbf{x}^{(K-1)}$ has already entered the quadratic convergence basin, reusing the transposed LU factors of $\mathbf{J}_{\mathbf{F}}(\mathbf{x}^{(K-1)})$ avoids an extra matrix factorization during the backward pass at the cost of a higher-order $O(\|\mathbf{x}^{(K)} - \mathbf{x}^{(K-1)}\|)$ perturbation.
\end{enumerate}
\end{remark}

\subsection{Numerical Experiments and SciML Benchmarking}
\label{subsec:experiments}

We evaluate the method across two 2-connected grid topologies and two physical friction regimes to test forward convergence, power-law generality, network scalability, and inverse parameter learning under both noiseless and noisy observations.

\subsubsection{Benchmark Topologies and Physical Regimes}
The benchmark suite comprises:
\begin{itemize}
\item \textbf{25-Branch Network ($4 \times 4$ Grid + Diagonal Trunk Chord):} $m = 16$ junction nodes ($\mathbf{A}_r \in \mathbb{R}^{15 \times 25}$), $n = 25$ branches, and $c = 10$ fundamental cycles ($\mathbf{C} \in \mathbb{R}^{10 \times 25}$), evaluated under both quadratic \textbf{Darcy--Weisbach ($\alpha = 2.0$)} and non-quadratic \textbf{Hazen--Williams ($\alpha = 1.852$)} friction laws.
\item \textbf{113-Branch Network ($8 \times 8$ Grid + Diagonal Trunk Chord):} $m = 64$ junction nodes ($\mathbf{A}_r \in \mathbb{R}^{63 \times 113}$), $n = 113$ branches, and $c = 50$ fundamental cycles ($\mathbf{C} \in \mathbb{R}^{50 \times 113}$), yielding a $113 \times 113$ coupled system.
\item \textbf{Multi-Scale Resistances:} Ground-truth resistances $\mathbf{s}^* \in \mathbb{R}_{++}^n$ span $[40.65, 5126.60]$ (condition ratio $p > 125$).
\end{itemize}

\subsubsection{Experiment 1: Forward Solver Convergence Across Power Laws and Scales}
We compare our dynamic-sign norm-linearization solver ($K=4$, on $F_i(\mathbf{x}) = u_i(\mathbf{x}) - v_i(\mathbf{x}) = 0$) against three root-finding baselines on the raw degree-$\alpha$ polynomial system $\widetilde{F}_i(\mathbf{x}) = u_i(\mathbf{x})^\alpha - v_i(\mathbf{x})^\alpha = 0$ \eqref{eq:raw_poly_res}: (i) \textbf{Undamped Newton ($K=4$)} without line search; (ii) \textbf{Standard Damped Newton} with Armijo backtracking line search ($\alpha_0 = 1.0$); and (iii) \textbf{Broyden's Quasi-Newton} method. On the 25-branch network ($1{,}000$ trials), initializations $\mathbf{x}^{(0)} = \lambda \mathbf{x}_{\text{lin}} \odot \exp(\boldsymbol{\xi})$ are generated around the oracle-free cold linear-resistive solution $\mathbf{x}_{\text{lin}}$ ($\alpha=1$) with global scale factor $\lambda \sim \mathcal{U}(2.0, 12.0)$ and log-uniform noise $\xi_k \sim \mathcal{U}(-0.22, 0.22)$. On the 113-branch network ($200$ trials), initializations apply the same scale dispersion and log-uniform perturbation inside the sign cone. Convergence is evaluated against $\|\mathbf{F}(\mathbf{x})\|_\infty \le 10^{-6}$ \eqref{eq:deq_operator}.

\begin{table}[h]
\centering
\footnotesize
\setlength{\tabcolsep}{4pt}
\caption{Forward-pass solver comparison across friction exponents ($\alpha=2.0$ and $\alpha=1.852$) and network scales ($n=25$ and $n=113$ branches).}
\label{tab:solver_comparison}
\resizebox{\textwidth}{!}{%
\begin{tabular}{llcccc}
\toprule
\textbf{Network / Regime} & \textbf{Solver Algorithm} & \textbf{Mean Iters} & \textbf{LS Evals} & \textbf{Non-Conv.} & \textbf{Time [ms]} \\
\midrule
25-Branch ($c=10$) & Undamped Newton ($K=4$, raw power-law) & 4.0 & 0.0 & $100.0\%$ & 0.59 \\
Darcy--Weisbach & Standard Damped Newton (Armijo) & 8.8 & 9.3 & $0.0\%$ & 7.15 \\
$\alpha = 2.0$ (Cold Start) & Broyden's Quasi-Newton & 41.4 & 0.0 & $60.1\%$ & 31.13 \\
& \textbf{Norm-Linearization ($K=4$, dynamic $\operatorname{sgn}$)} & \textbf{4.0} & \textbf{0.0} & \textbf{0.0\%} & \textbf{1.94} \\
\midrule
25-Branch ($c=10$) & Undamped Newton ($K=4$, raw power-law) & 4.0 & 0.0 & $100.0\%$ & 0.59 \\
Hazen--Williams & Standard Damped Newton (Armijo) & 8.0 & 8.3 & $0.0\%$ & 6.64 \\
$\alpha = 1.852$ (Cold Start) & Broyden's Quasi-Newton & 39.0 & 0.0 & $42.0\%$ & 30.33 \\
& \textbf{Norm-Linearization ($K=4$, dynamic $\operatorname{sgn}$)} & \textbf{4.0} & \textbf{0.0} & \textbf{0.0\%} & \textbf{1.96} \\
\midrule
113-Branch ($c=50$) & Undamped Newton ($K=4$, raw power-law) & 4.0 & 0.0 & $100.0\%$ & 2.76 \\
Darcy--Weisbach & Standard Damped Newton (Armijo) & 7.5 & 7.5 & $0.0\%$ & 28.74 \\
$\alpha = 2.0$ ($8\times 8$ Grid) & Broyden's Quasi-Newton & 41.2 & 0.0 & $55.0\%$ & 190.23 \\
& \textbf{Norm-Linearization ($K=4$, dynamic $\operatorname{sgn}$)} & \textbf{4.0} & \textbf{0.0} & \textbf{0.0\%} & \textbf{9.55} \\
\bottomrule
\end{tabular}
}
\end{table}

Across all three settings in Table~\ref{tab:solver_comparison}, undamped polynomial Newton fails on $100\%$ of trials at budget $K=4$ due to the non-vanishing affine shift $(\alpha-1)\widetilde{F}_i(\mathbf{x}^{(k)})$, whereas dynamic-sign Norm-Linearization converges below $10^{-6}$ on $100\%$ of trials in exactly $K=4$ unshifted linear solves.

\subsubsection{Experiment 2: Inverse Parameter Identification Under Prior Mismatch and Sensor Noise}
We next evaluate end-to-end inverse learning on the 25-branch network without oracle sign or warm-start leakage:
\begin{itemize}
\item \textbf{Prior-Centered Identifiable Parameterization ($\mathbf{z}^* \ne \mathbf{0}$):} We observe sparse flows on $8$ sensor branches ($\Omega_{\text{obs}} = \{7, 8, 12, 13, 14, 16, 17, 24\}$, $32\%$ coverage). To avoid centering the parameterization at the ground truth, resistances are parameterized relative to an imperfect nominal prior $\mathbf{s}_{\text{prior}} = \mathbf{s}^* \odot \exp(-\mathbf{B}_{\text{basis}}\mathbf{z}^*)$ ($30.63\%$ initial error on $\mathbf{s}$, where $\mathbf{z}^* \sim \mathcal{U}(-0.55, 0.55)^8$ and $\mathbf{B}_{\text{basis}}[0, :] = \mathbf{0}$):
\begin{equation}
\mathbf{s}(\mathbf{z}) = \mathbf{s}_{\text{prior}} \odot \exp(\mathbf{B}_{\text{basis}} \mathbf{z}), \quad \text{initialized at } \mathbf{z}^{(0)} = \mathbf{0}.
\end{equation}
Thus the Tikhonov regularizer $\gamma \|\mathbf{z}\|_2^2$ ($\gamma = 10^{-5}$) penalizes deviation from the imperfect prior $\mathbf{s}_{\text{prior}}$ rather than pulling $\mathbf{z}$ toward $\mathbf{z}^*$. Across $B=8$ demand snapshots, all implicit layers initialize their forward pass from an oracle-free cold linear-resistive start $\mathbf{x}_b^{(0)} = 2.5\,\mathbf{x}_{\text{lin}}(\mathbf{s}_{\text{prior}}, \mathbf{q}_b)$. We evaluate both \textbf{noiseless sensors ($0\%$ noise)} and \textbf{noisy sensors ($1\%$ multiplicative Gaussian noise)}.
\item \textbf{Compared Architectures (500 Adam Epochs):}
\begin{enumerate}
\item \textbf{Soft-Penalty MLP-PINN \cite{raissi2019}:} A 2-layer neural network $\mathbf{x}_{\mathbf{w}}(\mathbf{q}_r) = \mathbf{W}_2 \tanh(\mathbf{W}_1 \mathbf{q}_r + \mathbf{b}_1) + \mathbf{b}_2$ ($15 \to 32 \to 25$) maps nodal demands $\mathbf{q}_r$ to branch flows and is trained jointly with $\mathbf{z}$ under soft conservation penalties ($\lambda_{\text{phys}} = 10^3$, $\text{lr}=5\times 10^{-3}$).
\item \textbf{Raw Quadratic Newton-DEQ ($K=4$, on $u_i^2 - v_i^2 = 0$):} Runs $K=4$ undamped Newton steps on $\widetilde{\mathbf{F}}(\mathbf{x}, \mathbf{s}) = \mathbf{0}$ from the cold start $\mathbf{x}_b^{(0)}$ and differentiates via the quadratic IFT adjoint ($\text{lr}=5\times 10^{-3}$).
\item \textbf{Fixed-Budget FlowDEQ Layer ($K=4$, on $u_i - v_i = 0$):} Runs $K=4$ dynamic-sign norm-linearization steps from $\mathbf{x}_b^{(0)}$ and differentiates via the degree-one IFT adjoint \eqref{eq:vjp_closed} ($\text{lr}=3.5\times 10^{-2}$).
\end{enumerate}
\end{itemize}

\begin{table}[h]
\centering
\footnotesize
\setlength{\tabcolsep}{5pt}
\caption{Inverse parameter identification performance after 500 epochs from an imperfect prior ($30.63\%$ initial error) under $0\%$ and $1\%$ sensor noise.}
\label{tab:inverse_results}
\resizebox{\textwidth}{!}{%
\begin{tabular}{llccc}
\toprule
\textbf{Sensor Regime} & \textbf{Model Architecture} & \textbf{Rel. Error on $\mathbf{s}$} & \textbf{Max Phys. Violation} & \textbf{Time [s]} \\
\midrule
Noiseless & Soft-Penalty MLP-PINN ($\lambda_{\text{phys}} = 10^3$) \cite{raissi2019} & $27.40\%$ & $3.61 \cdot 10^{0}$ (Severe) & 0.12 \\
($0\%$ Noise) & Raw Quadratic Newton-DEQ ($K=4$, on $u_i^2 - v_i^2 = 0$) & $4.36\%$ & $9.44 \cdot 10^{-1}$ (Unconverged) & 1.68 \\
& \textbf{Fixed-Budget FlowDEQ Layer ($K=4$, on $u_i - v_i = 0$)} & \textbf{0.17\%} & $\mathbf{1.50 \cdot 10^{-7}}$ \textbf{($< 10^{-6}$ tol.)} & \textbf{1.35} \\
\midrule
Noisy Sensors & Soft-Penalty MLP-PINN ($\lambda_{\text{phys}} = 10^3$) \cite{raissi2019} & $27.40\%$ & $3.61 \cdot 10^{0}$ (Severe) & 0.13 \\
($1\%$ Gaussian) & Raw Quadratic Newton-DEQ ($K=4$, on $u_i^2 - v_i^2 = 0$) & $3.91\%$ & $1.00 \cdot 10^{0}$ (Unconverged) & 1.70 \\
& \textbf{Fixed-Budget FlowDEQ Layer ($K=4$, on $u_i - v_i = 0$)} & \textbf{1.12\%} & $\mathbf{1.60 \cdot 10^{-7}}$ \textbf{($< 10^{-6}$ tol.)} & \textbf{1.36} \\
\bottomrule
\end{tabular}
}
\end{table}

As shown in Table~\ref{tab:inverse_results}, the 2-layer MLP-PINN struggles to balance competing data and loop penalties ($27.40\%$ parameter error and $3.61$ conservation residual), while Raw Quadratic Newton-DEQ ($K=4$) retains an unconverged forward residual of $\sim 1.0$ that biases its IFT gradients (Remark~6.4.2). Dynamic-sign \texttt{FlowDEQ} ($K=4$) maintains physical residuals around $1.5 \cdot 10^{-7}$ throughout training, reducing parameter error from $30.63\%$ to $0.17\%$ on clean data and $1.12\%$ under $1\%$ sensor noise.

\section{Conclusions}

\begin{enumerate}
\item We derived the sharp universal upper bound $\eta_{\mathrm{univ}}(p) = \frac{(p-1)^2}{4p}$ for the relative variance functional alongside its exact finite-$n$ discrete counterpart, and proved via the weighted Kantorovich inequality (Theorem~\ref{thm:distortion_bridge}) that $\eta_{\mathrm{univ}}(p)$ is the exact worst-case single-step squared-norm distortion of the linearization operator when the target flow lies within a multiplicative condition ratio $p$ of the operating point.
\item Using Euler's homogeneous-function identity, we showed that the norm-linearization scheme is algebraically identical to Newton--Raphson applied to degree-one homogeneous norm differences $F_i(\mathbf{x}) = u_i(\mathbf{x}) - v_i(\mathbf{x}) = 0$, with an identically vanishing affine term, and proved via the row-equilibration identity $\mathbf{J}_{\text{loop}}(\mathbf{x}^*) = \mathbf{W}(\mathbf{x}^*)^{-(\alpha-1)}\mathbf{C}\mathbf{D}(\mathbf{x}^*)$ (Proposition~\ref{prop:nonsingularity}) that the coupled Jacobian is nonsingular at every regular equilibrium.
\item We proved that transforming raw degree-$d$ quadratic ($d=2$) or power-law ($d=\alpha>1$) loop equations $\widetilde{F}_i(\mathbf{x}) = u_i(\mathbf{x})^d - v_i(\mathbf{x})^d = 0$ into degree-one norm differences $F_i(\mathbf{x}) = u_i(\mathbf{x}) - v_i(\mathbf{x}) = 0$ acts as a nonlinear degree-reducing preconditioner (Proposition~\ref{prop:preconditioning}) that eliminates radial curvature ($\nabla^2 F_i(\mathbf{x})\mathbf{x} = \mathbf{0}$), cancels the affine residual shift $(d-1)\widetilde{F}_i(\mathbf{x}^{(k)})$, renders the loop Jacobian scale-invariant (degree-zero homogeneous), and row-equilibrates multi-scale circuits, verified across both $\alpha=2.0$ and $\alpha=1.852$ and up to $113$-branch networks.
\item Embedding this solver with dynamic sign partitioning into a Differentiable Physics Implicit Layer (\texttt{FlowDEQ})—backed by global variational uniqueness (Proposition~\ref{prop:global_uniqueness}) and closed-form IFT adjoint gradients (Theorem~\ref{thm:ift_adjoint})—enables fast, fixed-budget forward passes and accurate inverse parameter recovery from cold starts under both noiseless and noisy sensor observations.
\end{enumerate}

\section*{Code and Reproducibility}
The deterministic Python/PyTorch reference implementation (\texttt{main.py}) reproducing all theoretical verifications (Section~4 and Section~5), multi-regime and multi-scale forward solver benchmarks (Table~\ref{tab:solver_comparison}), and inverse parameter identification experiments under prior mismatch and sensor noise (Table~\ref{tab:inverse_results}) is provided alongside this manuscript.


\end{document}